\documentstyle[12pt,epsf]{article}

\newfont {\Bbb}{cmbx12}

\author{
D. V. Chistyakov \footnotemark
}

\title{
Fractal Geometry For Images Of Continuous Map
Of p-Adic Numbers And  p-Adic  Solenoids Into
Euclidean Spaces
}

\date{}
\begin{document}
\hyphenation{¥á-â¥-áâ-¢¥­-­®}
\thispagestyle{empty}
\newcommand{\bref}[1]{(\ref{#1})}

\newcommand{\Imm }{ \Upsilon_{s}^{(m)}}
\newcommand{\Immi }{ \Upsilon_{s}^{(\infty)}}
\newcommand{\Immn }[1]{ \Upsilon_{s}^{(#1)}}
\newcommand{\Immnn }[2]{ \Upsilon_{#2}^{(#1)}}

\newcommand{\wsol}{ \omega_{s}^{(m)}}
\newcommand{\wsoli}{ \omega_{s}^{(\infty)}}
\newcommand{\wsolip}{ \omega_{s/p}^{(\infty)}}
\newcommand{\Wsol}{ \Omega_{s,a}^{(m)}}
\newcommand{\Wsoli}{ \Omega_{s,a}^{(\infty)}}
\newcommand{\Wsolip}{ \Omega_{s/p,+0a}^{(\infty)}}
\newcommand{\Wsolipa}{ \Omega_{s/p,\ep   a}^{(\infty)}}
\newcommand{\idwsol}{ _{\otimes}\omega_{s}^{(m)}}
\newcommand{\facwsol}{  \mho_{s}^{(m)  }}
\newcommand{\wsoln }[1]{ \omega_{s}^{(#1)}}

\newcommand{\fp }[1]{ \{#1\}_p}
\newcommand{\ip }[1]{ [#1]_p}

\newcommand{\hin }[1]{ \chi_{n}^{(m)}(#1)}
\newcommand{\hi }[1]{ \chi_{#1}^{(m)}}
\newcommand{\hins }{ \tilde \chi_{n}^{(m)} }
\newcommand{\hinn}[2]{\chi_{n}^{(#1)}(#2)}
\newcommand{\hinsn}[1]{\tilde \chi_{n}^{(#1)} }
\newcommand {\norm}[2]{\mid #1 \mid_{#2}}
\newcommand {\np}[1]{\mid #1 \mid_{p}}
\newcommand {\npa}[1]{\mid #1 \mid_{p}^\alpha}
\newcommand {\abs}[1]{  \mid #1 \mid  }
\newcommand {\half}{ \frac{1}{2} }
\newcommand{\EL}{\stackrel{L}\simeq}
\newcommand{\Ll}{\stackrel{L}\preceq}
\newcommand{\Lg}{\stackrel{L}\succeq}
\newcommand{\proof}{{\bf Proof.~}}
\newcommand{\prooflem}{{\bf Proof of the lemma.~}}
\newcommand{\cond}{{\bf Convention~}}

\newcommand{\Endproof}{$\Box$}
\newcommand{\remark}{³àìå÷àíèå: }
\newcommand{\Ran }{{\rm Ran}}
\newcommand{\RE }[1]{{\rm Re}\left ( #1 \right )}
\newcommand{\IM }[1]{{\rm Im}\left ( #1 \right )}
\newcommand{\Rs}{{\bf R}}
\newcommand{\Su}{{\bf S_1}}
\newcommand{\sC}{{\bf C}}
\newcommand{\Zs}{{\bf Z}}
\newcommand{\Qs}{{\bf Q}}
\newcommand{\Cs}{{\bf \C}}
\newcommand{\Ns}{{\bf N}}
\newcommand{\eNs}{\bar {\bf N}}
\newcommand{\Tp}{{\bf T}_{p}}
\newcommand{\Qp}{{\bf Q}_{p}}
\newcommand{\Zp}{{\bf Z}_{p}}
\newcommand{\Tn}[1]{{\bf T}_{#1}}
\newcommand{\Qn}[1]{{\bf Q}_{#1}}
\newcommand{\Zn}[1]{{\bf  Z}_{#1} }
\newcommand{\FRp}[1]{ \Rs_{+}^{#1} }
\newcommand{\pic}[2]{ fig. \ref{#1}.#2 }
\newcommand{\pf}[1]{ fig. \ref{fpi}.#1 }
\newcommand{\ps}[1]{ fig. \ref{spi}.#1 }

\newcommand{\UC}{\rm U_1}
\newcommand{\diam}{\rm diam}
\newcommand{\la}{\lambda}
\newcommand{\ep}{\varepsilon}
\newcommand{\tilF}{\widetilde{F}}
\newcommand{\tilP}{\widetilde{P}}
\newcommand{\tilG}{\tilde G}
\newcommand{\tilC}{\widetilde{C}}
\newcommand{\tilT}{\widetilde{T}}
\newcommand{\tila}{{\tilde a}}
\newcommand{\tilf}{{\tilde f}}
\newcommand{\tilY}{\widetilde{Y}}
\newcommand{\tilZ}{\widetilde{Z}}
\newcommand{\tilO}{\widetilde{O}}
\newcommand{\tilPhi}{\widetilde{\Phi}}
\newcommand{\myind}{\mbox{\rm Ind\,}}
\newcommand{\FI}{$\,\mbox{L.L.}\!$}
\newcommand{\FIp}{\mbox{L.L.$p$-L}}
\newcommand{\tilR}{\widetilde{R}}
\newcommand{\calR}{{\cal R}}

\newtheorem{mydef}{Definition}
\newtheorem{mylemma}{Lemma}
\newtheorem{mytheorem}{Theorem}

\maketitle
\begin{abstract}
\noindent
\footnotetext{Kazan' State University. E-mail: Dmitry.Chistyakov@ksu.ru}


Explicit formulas are obtained for a family of continuous mappings of p-adic
numbers $\Qp$  and solenoids $\Tp$  into the complex plane $\sC$ and the space
~$\Rs ^{3}$  , respectively. Accordingly, this family includes the mappings for
which the Cantor set and the Sierpinski  triangle  are images of the unit balls
in
$\Qn{2}$  and  $\Qn{3 }$. In each of the families, the subset of the embeddings
is found. For these embeddings, the Hausdorff dimensions are calculated and it
is shown that the fractal measure on the image of $\Qp$ coincides with the Haar
measure on $\Qp$. It is proved that under certain conditions, the image of the
$p$-adic solenoid is an invariant set of fractional dimension for a dynamic
system. Computer drawings of some fractal images are presented.

\end{abstract}
\section
{Introduction}\label{Introduction}

The hierarchical structure of $p$-adic numbers and fractals, as well as the
symmetries of self-similar fractals, point to a close relationship between
these objects that has been repeatedly noted  \cite{VVZ,Pit,dif1}. A clear
example of this kind is the homeomorphism of the Cantor set onto the ring
$\Zn{2} $ ~ \cite {VVZ,Zelenov,Pit}. However, as far as we know, explicit
formulas for the embeddings of various subsets of $\Qp$ into Euclidean spaces
and the fractal properties of the corresponding images have received little
attention. At the same time, the topology of these objects as strange
attractors  \cite{LL} is similar to that of the $p$-adic solenoids $\Tp$
and therefore, the construction of the embeddings of $\Tp$ , say, into the
three-dimensional Euclidean space  $\Rs^{3}$  is also of interest.

In this paper we construct the continuous mappings  $\Imm : \Qp \mapsto \sC$,
depending on
the parameters $s,a \in \sC$ and the number m, which can be $a$ positive
integer or $\infty$. It turns out that $\Immnn{0}{1/3}(\Zn{2})$ is a Cantor
set and $\Immnn{0}{1/2}(\Zn{3})$  is a Sierpiriski carpet. It is shown that
for some $s_0>0$, the condition $\abs{s} < s_0$
( ¨ $\abs{a} < (1-\abs{s})^{-1}$) determines the sets of mappings
$\Imm $ and $\Wsol$  that are the embeddings.

These sets of embeddings possess the following property. Given $s$ and $a$,
there exist additively invariant metrics in $\Qp$ and $\Tp$ such that the
mappings $\Imm $ and  $\Wsol $ a preserve the Hausdorff dimensions of arbitrary
subsets. This property proves to be convenient when studying the Hausdorff
measures of the sets.
This property proves to be convenient when studying the Hausdorff measures of
the sets
$\Imm(\Qp)  $ and  $\Wsol(\Tp)$  because the Hausdorff measures in the spaces
$\Qp$ and $\Tp$, with the corresponding metrics. are simply the Haar measures.
Moreover, it turns out that the Haar measure of any set in $\Qp$ coincides
with the fractal measure of its image (see formula\bref{frmes}).
This property and the fact that $\Immi$  is a series of continuous additive
characters in $\Qp$ make it possible to apply group-theory methods when
calculating the integrals with respect to the fractal measure. This can serve
as an effective means for solving some problems of quantum mechanics and
diffusion or diffraction on such fractals  \cite{dif1,dif2,FFT}.

We prove that for an embedding $\Wsoli $ there is a dynamic system
\footnote{By a dynamic system in $\Rs^n$  we mean an autonomous system of $n$
first-order equations that satisfies the conditions of the existence and
uniqueness theorem}
such that $\Wsol(\Tp)$  is an invariant set of this system and any integral
trajectory lying in $\Wsol(\Tp)$ densely winds around this set. In addition,
this mappings
$\Immi $ and $\Wsoli $, are inter-related. In the present paper, we construct
an injective homomorphism $j$ of the additive group $\Qp$ into the group $\Tp$
such that $j(\Qp)$ is dense in
$\Tp $ and prove the existence of a Lipschitz mapping $J: \sC \mapsto \Rs^3 $
commuting with $j$ (see formula \bref{slice}) whose restriction to $\Imm(\Qp)$
is a local isometry.

\section
{
Hausdorff measures (basic definitions)
}\label{hmess}
Here we present some relevant material concerning Hausdorff measures in a form
that is suitable for the subsequent presentation. For an arbitrary metric space
$(M,\rho)$
, we define the $\delta$-dimensional outer Hausdorff measure$h^{\delta}$ by
setting  $ \forall A \subset M$,
\begin{eqnarray}
h^{\delta }(A)=\lim_{\ep \to +0}  h_{\ep }^{\delta}(A)
=\lim_{\ep \to +0}  \inf
\{
\sum_{i=1} \diam(S_{i})^{\delta} :
\bigcup_{i=1}^{\infty} S_{i} \supseteq A , \diam(S_{i}) \leq {\ep}
\},
\label{mhe}
\end{eqnarray}
where $\delta$ is a fixed positive number and $ \diam (B) \equiv  \sup
\{ \rho (x,y)  : x,y \in B \}$ . Then $ h^{\delta}$ is a measure countably
additive and regular in the sense of Borel \cite{FED}. By definition, the
Hausdorff dimension of a subset $ \forall A \subset M$  is the number
\cite{bill}
\begin{eqnarray}
D_{h}(A)=\inf \{ \delta :  h^{\delta}(A)= 0 \} =\sup \{ \delta :
h^{\delta}(A)=  \infty \}.
\end{eqnarray}

In particular, it follows that if $  h^{\delta}(A) >0 $ and the measure
$h^\delta_{\mit A}(\cdot)\equiv h^\delta(\cdot \cap A)$~is $\sigma$-finite,
then $\delta =  D_{h}(A) $. By the local Hausdorff dimension at a point
$ x \in A$, we mean the number $D_{h}^{L}(x)=\inf (D_h(A\cap U_x)),$
where $\inf$ extends over all open neighborhoods  $U_x$ of the point  $x$.

\cond. Irrespective of the nature of the set $X$, for arbitrary non-negative
real-valued functions $F$
and $G$ on $X$, we write $ F \Ll G$  or $  \forall x \in X~ F(x) \Ll G(x) $
wherever there exists $C > 0$ such that  $  \forall x \in X~ F(x) \leq CG(x) $.
If the relations $ F \Ll G$   and  $G \Ll F$ hold simultaneously, then
$ F \EL G$ and we say that F and G
are equivalent
\footnote{Actually, the symbols $\Ll$ and $\EL$
define order and equivalence relations, respectively.}
.
When interpreting constants as trivial functions, we write $c \EL 1$ instead
of $0< c < \infty$.
For any  $F,G$, and $H \Ll F$. the following elementary relations hold:
$ \forall a,b,\alpha >0 $
\begin{eqnarray}\label{LEQ}
\min (F,H) \EL H, \max (F,H) \EL F,
\nonumber \\
~aF+bG \EL \max(F,G) \EL  \max(F,G \pm H)  \EL  ( F^{\alpha}+G^{\alpha}
)^{\frac{1}{\alpha}}.
\end{eqnarray}

Let $\Phi : M \mapsto N$  be a mapping of the metric spaces $(M,\rho)$ and
$(N,d)$
and let  $d^{\Phi}\equiv d(\Phi(\cdot),\Phi(\cdot))$
Then $\Phi$  is a Lipschitz mapping if $d^{\Phi} \Ll \rho$. In this case.
$\Phi$ is called an {\it L-contraction}.
If  $d^{\Phi} \Ll \rho$. then   $\Phi$  is, a Lipschitz embedding and the
restriction of $\Phi^{-1}$  to $\Ran(\Phi)$
is also a Lipschitz mapping: such a mapping is called a
{\it Lipschitz isometry}
or {\it L-isometry}.

\begin{mydef}\label{luft}
For given pseudometrics ~$\rho_1$~and ~$\rho_2$ on a set $M$, the quantity
\begin{eqnarray}\label{kappa}
\kappa_{1,2}= \sup_{ x,y \in M} \left ( \frac{ \abs{ \rho_1(x,y) -
\rho_2(x,y)  }}{ \rho_1(x,y)  + \rho_2(x,y) } \right )
\end{eqnarray}
( where we set  $ 0/ 0 =0$) is called the divergence of ~$\rho_1$~and
~$\rho_2$.
\end{mydef}
Obviously,  $\kappa_{1,2} \leq 1$ and, moreover, it can be shown that {\it
$\kappa_{1,2} < 1$~
if and only if  $ \rho_1 \EL  \rho_2$.}
It is easy to prove that the following assertion holds.
\begin{mytheorem} \label{eqmet}
Let $ \rho_1$ and $\rho_2$- be metrics on $M$, let $\kappa_{1,2} $ be their
divergence,
and $ h^{\delta}_1,D_1,D_1^{L}$ and $ h^{\delta}_2,D_2,D_2^{L} $   denote
the corresponding
$\delta$-dimensional measures and global and local Hausdorff dimensions.
Then the inequality
\begin{eqnarray}
\left (\frac{1-\kappa_{1,2}}{1+\kappa_{1,2}} \right)^{\delta}  h^{\delta}_1(F)
\leq
h^{\delta}_2(F)
\leq
\left (\frac{1+\kappa_{1,2}}{1-\kappa_{1,2}} \right)^{\delta}  h^{\delta}_1(F),
\end{eqnarray}
holds, whence follows that if $\rho_1 \EL \rho_2$, then $h^{\delta}_1 \EL
h^{\delta}_2$
and we have $\forall F \subset M ~,
D_{h,1}(F) =D_{h,2}(F) $ and $ D^{L}_{h,1}(x) =D^{L}_{h,2}(x)~ \forall x
\in F $.
\end{mytheorem}

In view of  \ref{eqmet}, Theorem 1 and Theorem 2.10.45 of  \cite{FED} imply
the next theorem.
\begin{mytheorem} \label{RxM}
Let $(M,\rho)$ be a metric space. Introduce a metric $d$ in $ \Rs \times M $
by the formula
$ d((x,a),(y,b))=\max ( \abs{x-y}, \rho (a,b) ),~  \forall x,y \in
\Rs ~¨~\forall a,b \in M$.
Then for every Borel set $ A \subset \Rs$ and every $\forall B \subset M$ ~
with~$h^{\delta}(B)< \infty $~  , we have
\begin{displaymath}
h^{ \delta +1  }(A \times B) \EL h^{1}(A) h^{ \delta }(B ).
\end{displaymath}
\end{mytheorem}

\section
{
Hausdorff measures in $\Qp$~ and~$\Tp$
}\label{padic}

In this section,  $\Qp$ and $\Tp$ are regarded as completions of $\Qs$
with respect to the corresponding additively invariant metrics.
Each element $x$ of the $p$-adic number field $\Qp$ is uniquely representable
as a formal series  \cite{Kob},
\begin{eqnarray}\label{power}
x=\sum_{n=v}^{\infty}a_{n} p^{n}
= \sum_{n=v}^{-1}a_{n} p^{n}+\sum_{n=0}^{\infty}a_{n} p^{n}
\end{eqnarray}
with coefficients $a_{n} \in \{0,1,...,p-1\}$ where $v<\infty$ and $p$ is a
fixed prime number.
\footnote
{
Actually, the part of $p$ can be played by any positive integer because we
do not use the existence of inverse elements the ring $\Qp$ anywhere.
}
This series absolutely converges in the $p$-adic norm defined $\forall x$
by the relation $ \Vert x \Vert = p^{- \alpha v(x)}$
for some $\alpha >0$ and $v(x)=v$ is called the logarithmic norm of $x$.
The first sum on the right-hand side of \bref{power} is denoted as
$\fp{ x }$ and is the fractional part of $x$. and the other is denoted as
$\ip{x}$ and is the integral part of $x$. In this case,
$\fp{x} \in \Qs \cap [0,1) $ ~ and $\fp{x} \in \Qs \cap [0,1) $ ~, where
~$ \Zp=\{ x \in \Qp:  \Vert x \Vert \leq 1  \} $~ is the ring of integer
$p$-adic numbers. Any number $q \in \Qs$ can be expanded uniquely as a
series \bref{power} and $\Qp$ is the completion of $\Qs$ \cite{Kob}.
The norm  $\Vert \cdot \Vert$  possesses the property of ultrametricity,
$ \forall x,y \in  \Qp$, we have
\begin{eqnarray}\label{ultra}
\Vert  x-y  \Vert  \leq  \max( \Vert x  \Vert ,  \Vert y  \Vert ).
\end{eqnarray}

The norm with~ $ \alpha =1$ ~is denoted by~ $\np{\cdot  }$ and the others
are denoted simply as ~ $\np{\cdot  }^{\alpha }$.  All of these norms
depending on a are topologically equivalent  \cite{Kob}; however, the
corresponding Hausdorff measures are different (see below).
The standard Haar measure$\chi$ in  $\Qp$ is chosen such that \cite{Kob}
\begin{eqnarray}\label{nhaar}
\chi  (\Zp) =  \int  \limits_{ \Zp} d \chi = 1.
\end{eqnarray}
\begin{mytheorem}
The Hausdorff measure $h^{1/\alpha}$ on $(\Qp,\np{\cdot}^{\alpha})$ is
the standard regular Haar measure $~\forall \alpha >0$ and, consequently,
the local and global Hausdorff dimensions of
$(\Qp,\np{\cdot}^\alpha)$ coincide and are equal to $1/ \alpha$.
\end{mytheorem}
\proof
By construction, $h^{\delta}$  is an invariant measure $\forall \delta > 0$~
and, therefore, the uniqueness of
$\chi$ implies that it suffices to show that, for instance,
$h^{1/\alpha}(\Zp) =\chi(\Zp)=1$. Because $\Zp$ consists of
exactly $p^{M}$ disjoint balls of diameter  $p^{-M \alpha }~\forall M \in \Ns$,
we have $h^{1/\alpha}(\Zp) \leq 1.$
We now show that $h^{1/\alpha}(\Zp) \geq 1.$  Indeed, by the semiadditivity of
$\chi ~ \forall A \subset \Qp, \forall  \ep>0$ , the inequality
\begin{eqnarray}   \label{h000}
\chi(A) \leq \inf
\{
\sum_{i=1} \chi(U_{i}) :
\bigcup_{i=1}^{\infty} U_{i} \supseteq A , \diam(U_{i}) \leq {\ep}
\},
\end{eqnarray}
holds, where $ U_{i}$, are open balls in $\Qp$. On the other hand, it follows
from the properties of norm ($\ref{ultra}$) that any subset  $ B \subset
\Qp$~ á ~$\diam(B) = r $ ~
is contained in the open (and, simultaneously, closed) ball $U=\{ x :
\np{x-x_{0}}^{\alpha} \leq r , x_{0} \in B \}$
with ~ $\diam(U)=r $ and this, together with the relation $\chi(U)=
\diam^{1/\alpha}(U)$,
which is valid for any ball $U$, implies that $h_{\ep}^{1/\alpha}$  coincides
with the right-hand side of inequality
(\ref{h000}) and, consequently, we have $h^{1/\alpha}(\Zp)
\geq \chi(\Zp)$.\Endproof

Let us consider $\Rs  \times \Zp $  as a direct product of additive groups
and introduce a metric in this group by fixing some ~$\alpha > 0$~ and
setting ~$ \forall a,b \in \Rs ~ and ~ \forall x,y \in \Zp $
\begin{eqnarray}   \label{rxz}
\hat \rho_\alpha ((a,x),(b,y)) = \max ( \abs{a-b},\np{x-y}^\alpha ).
\end{eqnarray}
It is clear that $  \hat \rho_\alpha$  is an invariant metric and the
topology generated by it coincides with that of the direct, product of
groups.
\begin{mydef}\label{soldef}
Let  ~$B$ denote the subgroup $  \{ (n, - n):   n \in \Zs  \}$
of the group $\Rs  \times \Zp$. Then the quotient group  $ \Tp =  (\Rs
\times \Zp )/B$ is called
a $p$-adic solenoid
\footnote{
This definition differs from the one in  \cite{HR}, where $B =\{(n,n): n\in
\Zs \}$.
However, since the mapping $ x \mapsto -x $ is automorphism of the additive
group $\Zp$. the corresponding quotient groups are isomorphic.
}.
\end{mydef}
It can be shown  \cite{HR} that $\Tp$ is a connected compact Abelian group.
We can define an invariant metric $\rho_\alpha $~on $ \Tp $ ~ that is
compatible with the topology as the quotient metric according to the
standard scheme\cite{HR}. namely, $ \forall f,g \in \Tp~ $, we set
\begin{eqnarray}   \label{rt}
\rho_\alpha (f,g) = \inf \{ \hat \rho_\alpha ((a,x),(b,y)) :  (a,x) \in f ,
(b,y) \in g\}.
\end{eqnarray}
The construction below gives a concrete realization of ~$\Tp$.
\begin{mytheorem} \label{cansol}
Consider the product $[0,1) \times \Zp$ and define addition in it by the
following rule:
$\forall f=(\xi ,x),g=(\eta ,y) \in [0,1) \times \Zp $
\begin{eqnarray} \label{addsol}
f+g =( \xi + \eta- [\xi + \eta], x+y + [\xi + \eta]),
\end{eqnarray}
where $[\xi+\eta]$ is the  integral part of the real number $\xi+\eta$.
For a fixed $\alpha >0$, we define a metric $\rho_{\alpha}$ by setting
\begin{eqnarray}\label{rosol}
\rho_{\alpha}(f,g)=\min( \ell_{\alpha}(f-g), \ell_{\alpha}(g-f)),
\end{eqnarray}
where $\ell_{\alpha}(f)=\max( \xi ,\np{x}^\alpha) $.
Then the resulting  Abelian group with the topology induced by the metric is
algebraically and isometrically isomorphic $\Tp$ with metric  \bref{rt}.
\end{mytheorem}
\proof
The algebraic isomorphism is established in practically the same way as in
the proof of theorem 10.15  \cite{HR}.
It follows from \bref{rxz}, \bref{rt} and the definition of subgroup $B$ that
the metric on  ~$\Tp$ satisfies the relation
\begin{displaymath}
\rho_\alpha (f,0) = \inf \{ \max ( \abs{\xi -n },\np{x+n}^\alpha ).  : n \in
\Zs \},
\end{displaymath}
However, because $\forall x \in \Zp $, the inequality $\np{x+n}^\alpha \leq 1$
holds,
the infimum in the above formula can be extended over the set $n=0,1$.
Taking into account that $ ( -f)=( 1-\xi,-x-1) $ for $\xi \ne 0$,
we can easily show that the metrics do, in fact, coincide.
\Endproof

Using Theorem \bref{cansol}  and the facts that $\forall x \in \Zp
~ \rho_{\alpha}((0,x),0) = \np{x}^{\alpha}~$
and  $\fp{ x }=0 $, and, also, that the set $\{ \fp{x} : x \in \Qp \} $
is dense in  $[0,1)$  (in the usual topology), we can prove the following
theorem.
\begin{mytheorem}\label{JQT}
The mapping $ j: \Qp \mapsto \Tp $, associating the element  $(\fp{x},\ip{x})
\in \Tp$
with each $x \in \Qp$,
is an injective homomorphism of the additive group $\Qp$ into $\Tp$
and is also a local
isometry from $(\Qp,\np{\cdot}^{\alpha})$ into ~  $(\Tp,\rho^{\alpha})~\forall
\alpha >0$.
Furthermore,  $\Ran(\Qs) \subset \Ran(\Qp)$ is dense in $\Tp$.
Because $\Tp$ is complete (in view of the compactness \cite{Kel}),
$\Tp$  is the completion of $\Qs$ with respect to the metric $\rho^{\alpha}$
\footnote
{
Note that although $j$ is not an embedding ( $\Ran(\Qp)$ ¯«®â­® ¢ $\Tp$ ),
the restriction of $j$ to $p^{M} \Zp$ is an (L-) isometry  $\forall M \in
\Ns(\Zs)$.}.
\end{mytheorem}

The $p$-adic solenoid $\Tp$ is a compact Abelian group \cite{HR}
and, therefore, there is a finite Haar measure $\chi$ on it that is unique and
invariant.
Let us show that the Hausdorff measure $h_{\delta} $  on$(\Tp,\rho_{\alpha})$
for $\delta = \alpha^{-1}+1$ coincides
(to within a finite nonzero multiplier) with $\chi$ on all Borel subsets of
$\Tp$.

Since $\rho_{\alpha}$  is an invariant metric, it suffices
to prove that
$U_{1/2}=\{y \in \Rs \times \Zp : \rho_{\alpha}(y,0)
< 1/2 \}~ h_{\delta}(\phi(U_{1/2})) \EL 1$,
where ($\phi$ is the canonical projection of $ \Rs  \times \Zp$ onto $\Tp$.

However, because  $\forall a,b \in B \subset \Rs  \times \Zp $ ~
we have $  \hat \rho_\alpha (a,b) \geq 1 $ for $a \ne b$,
the restriction of $\phi$ to $U_{1/2}$ is an isometry.
Therefore, applying the relation $U_{1/2}=(-\frac{1}{2},\frac{1}{2})
\times p^{n_{\alpha}}\Zp$
for some $n_{\alpha}>0$,
we conclude, in view of Theorem (\ref{RxM}), that
\begin{displaymath}
h_{\delta }
(\phi(U_{1/2})) =
\hat h_{\delta }( (-1/2,1/2) \times \{ x \in \Zp : \np{x} \leq p^{-1} \} )
\EL p^{-n_{\alpha}} \tilde h_{ 1/ \alpha}( \Zp).
\end{displaymath}

Here  $\hat h_{\delta }$ and $\tilde h_{\delta }$
denote the Hausdorff measures in $ \Rs  \times \Zp$ and $\Tp$, respectively.
Thus, as in the case of $\Qp$, the {\it local and global Hausdorff dimensions
of $(\Tp,\rho_{\alpha})$ coincide and are equal to $ \alpha^{-1}+1$.}

\section
{Continuous mappings of  $\Qp$  into $\sC$
}\label{immers}
{\Bbb C}
Given  $\forall n \in \Zs $ and ~$\forall m \in \eNs \equiv  \Ns \cup
\{\infty\} $~
, we define the complex-valued functions$\hin{\cdot }$ on $\Qp$ by the formula
\begin{eqnarray} \label{hina}
\hin{x} =\exp (\frac{ i2 \pi }{p} \sum_{k=0}^{m} x_{n-k} p^{-k} ),
\end{eqnarray}
where $x_{n}$ is the $n$th coefficient in the expansion of $x$ into series
(\ref{power}).
It is easy to show that $\hin{\cdot }$ are
continuous. Note that $ \chi^{(\infty)}_{n}(\cdot) $
coincides with the continuous additive character
$\chi_{\frac{1}{p^{n+1}}}(\cdot)$  on $\Qp$
\cite{GG,HR}.
\begin{mydef}\label{mapQC}
For every $ s \in \UC \equiv \{ z\in \sC : \abs{z}<1 \}$
and  ~$\forall m \in \eNs $~,
we define a continuous mapping ~$\forall m \in \eNs $~ by setting
$ \Imm : \Qp \mapsto \sC  $
\begin{eqnarray} \label{impow}
\Imm (x)  =\frac{1- s^{v(x)} }{1-s} +\sum_{n=v(x)}^{\infty} s^{n}
\hin{x}=[\Imm](x)+\{ \Imm \}(x),~~~\forall x \in \Qp.
\end{eqnarray}
Here $[\Imm](x)=\sum_{n=0}^{\infty} s^{n} \hin{x}$ is the "integral part"
$\Imm(x)$
and $\{ \Imm \}=\Imm(x)-[\Imm](x)$  is the "fractional part" of $\Imm(x)$.
\end{mydef}

Obviously, the mapping $\Imm$ is well defined and for any fixed $x \in \Qp$,
formula \bref{impow}
defines a function of $s$ which is holomorphic in the circle $\UC$.
In addition, taking into account that the coefficients
of series  \bref{power} are zero (periodic \cite{Kob})
beginning with some $\forall q \in \Ns (\Qs) $ , we can show that the following
interesting
property holds:  $\Immnn{m}{(\cdot)}(q)-\Immnn{m}{(\cdot)}(0)$~
is a polynomial (a rational function) for $m < \infty  $
and  $\varrho(\cdot) \cdot \Immnn{\infty}{(\cdot)}(q)$ is an entire function
(a meromorphic function, i.e., a ratio of entire functions) for $m=\infty$,
where
$$\varrho(s) = \prod_{k=0}^{\infty}(1-p^{-k}s)$$.

In the case $m<\infty$ , the proof of this assertion follows from the
periodicity of $\hin{q}$
with respect to $n$ for  $n > n_q+m$. If  $m=\infty$,
we have
$$ ~ \chi^{(\infty)}_{n}(q)= \exp(i2\pi \fp{ q/p^{n+1} })=\exp(i2\pi q/p^{n+1})
\exp(-i2\pi \ip{ q/p^{n+1} } )$$,
$\forall q \in \Qs~$ and, therefore, it is possible to construct the analytic
continuation with the aid of the Taylor expansion of $\exp(i2\pi q/p^{n+1})$,
using the periodicity of $\exp(-i2\pi \ip{ q/p^{n+1} } )$
with respect to $n$ for  $n > n_q$.
\footnote
{The analytic properties of $\Immnn{m}{(\cdot)}(x)$ for $x \in \Qp \backslash
\Qs$
are unknown to the author. However, if in expansion  \bref{power}, say, we
set  $x_n = \sum_{k=0}^{\infty} \delta_{n,2^k}$
for $x \in \Qn{2}$, then  $\Immnn{0}{s}(x)=-2 \sum_{k=0}^{\infty}s^{2^k}$
and it can be shown \cite{Shabat} that $\UC$ is a holomorphy domain of
$\Immnn{0}{(\cdot)}(x)$
. }.
It is easy to prove that $\forall x \in \Qp$, we have the relations
\begin{eqnarray}
\{ \Imm \}(x) = \{ \Imm \}(
\fp{x}),~\Imm(\ip{x})=[\Imm](\ip{x}),~[\Imm](p^{-m}x)=[\Imm](p^{-m}\ip{x})
,\label{fract}
\end{eqnarray}
The following basic property of the mapping {\it (the scaling)} $\Imm$ also
holds:
$\Imm$ -{\it scaling}:
\begin{eqnarray}
\forall x \in \Qp~~ \Imm (px) = s\Imm(x)+1 =p^{-\frac{1}{D_{s}}} e^{i \arg(s)
} \Imm(x) +1 \label{scale},
\end{eqnarray}
where $ D_{s} = - {\rm log}_{p}^{-1}(\abs{s})$ is called the scaling dimension
of  $ \Imm $ .
This relation follows from the fact that
$ \chi^{(m)}_{n+k}(p^k x) =\chi^{(m)}_{n}(x) , \forall k \in
\Zs \label{xiprop1}$.
Furthermore, by applying \bref{fract} , \bref{scale}, we can show that the
set {\it $\Ran \Imm$   is self-similar
$\Ran \Imm$}.More precisely, let $ B_l^{n} \equiv \Imm( \{x \in \Qp : \np{x-l}
\leq p^{-n}\} )$, then
\begin{eqnarray}\label{selfsim}
B_l^{n} = \bigcup_{\bar l =0}^{p^m-1}
\{ z_{l,\bar l}^n +e^{i \arg(s)mn } p^{-\frac{mn}{D_{s}}} B_{\bar l}^0  \},
\end{eqnarray}
$\forall n \in \Zs, l \in \Qp $.
Here $ z_{l,\bar l}^n$ are shifts of$\sC$ depending on $l,\bar l$ and $n$.
Thus, every $B_l^{n}$ can be obtained
from $p^m$ sets $ B_{\bar l}^0 $ by means of continuous motions of the plane
$\sC$~ (shifts and rotations) and a scaling transformation.

Now, we find the conditions under which the mappings $\Imm$   become
embeddings.
To this end, we define the number
\begin{eqnarray} \label{da}
\Delta^{(m)}_s= \inf  \{ \abs{  \Imm(x) - \Imm(y) } :\forall x,y \in \Qp :
\np{x-y}=1\}.
\end{eqnarray}
Using the fact that $~\abs{ \chi^{(m)}_{0}(x) -  \chi^{(m)}_{0}(y)  } \geq 2
\sin( \pi /p)$
for$\np{ x - y} = 1$. we can derive the inequality
\begin{eqnarray}\label{ds0}
\frac{ \Delta^{(m)}_s }{2} \geq \sin \left ( \frac {\pi}{p} \right ) -
\frac{\abs{s}}{1-\abs{s}}.
\end{eqnarray}
It follows that $  \Delta^{(m)}_s>0$  for
$$
\abs{s} < s_0 =\frac{\sin(\pi/p)}{1+ \sin(\pi/p)}
$$
\begin{mytheorem}\label{LipImm}
Let $s$ and $m$ be such that $ \Delta^{(m)}_s >0$,
then $\Imm$ is an L-isometry from $(\Qp,\np{\cdot}^{D_{s}^{-1}})$
into $(\sC,\abs{\cdot})$
and, therefore, $\Imm$   is an embedding
\footnote{
Using the compactness of $\Zp$, property (\ref{scale}) and the
completeness of $\Qp$~ and ~$\sC$, we can show that if $\Imm$ is an injective
mapping, then it is an embedding
.}.
\end{mytheorem}
\proof
The proof can be derived from the inequality below, which is a consequence of
\bref{scale},
and definition  \bref{hina}, namely, $\forall x,y \in \Qp$ we have
\begin{eqnarray}\label{AA}
\Delta^{(m)}_s \abs{s}^{v(x-y)} \leq \abs{  \Imm(x) - \Imm(y) } \leq
\frac{2}{1- \abs{s}} \abs{s}^{v(x-y)} \mbox{.\Endproof}
\end{eqnarray}
It follows from Theorems   \ref{eqmet} and \ref{LipImm} that
{\it if  $ \Delta^{(m)}_s >0$, then the local and global Hausdorff dimensions
of $\Imm(\Qp)$ are equal to $D_s$.} These theorems also imply that
$ h^{D_s}(\Imm(\Zp)) \EL 1$.
Therefore, we can introduce a {\it fractal measure} $\mu_f$ on $\Qp$, by the
formula
\begin{eqnarray}\label{frmes}
\mu_f (B) =\frac{  h^{\delta}(B \cap \Imm(\Qp)) }{h^{\delta}(
\Imm(\Zp))}~~\forall B \subset \sC.
\end{eqnarray}
The restriction of this measure to $\Imm(\Zp)$ coincides with
the multifractal measure $\mu_f^0$ \cite{Feder} of the set $\Imm(\Zp)$.
Moreover, the assertion below holds.
\begin{mytheorem}\label{mueqxi}
Let $\chi(\cdot)$ be the standard Haar measure in $\Qp$. In this case, if
$ \Delta^{(m)}_s  > 0$, then
\begin{eqnarray}
\mu_f(\cdot)=\chi ( (\Imm)^{-1} (\cdot)) \label{muxi}.
\end{eqnarray}
\end{mytheorem}
\proof
Because  $\Imm$ is an embedding, the clusters $B_l^{n}$ defined in
(\ref{selfsim}) do not intersect.
Therefore. if $ m < \infty$, then the self-similarity of the set $\Imm(\Qp)$
implies that  $ h^{D_s}(\Imm(B))= h^{D_s}(\Imm(\Zp ))\chi(B)$ for all open
sets B. The proof of the theorem
in the case m = $\infty$ follows from Theorem  \bref{eqmet} and the lemma
below.
\begin{mylemma}
Let $\rho_m$ be pseudometrics on  $\Qp$ denned by the formula
$\rho_m(x,y)=\abs{\Imm(x)-\Imm(y)}$~$ \forall x,y \in \Qp$
and let $\kappa_{m,\infty}$ be the divergence of $\rho_m$~and~$\rho_{\infty}$.
Then
\begin{eqnarray*}
\kappa_{m,\infty} < \frac{4 \pi}{(1-\abs{s})( \Delta^{(m)}_s+
\Delta^{(\infty)}_s)}p^{-m} .
\end{eqnarray*}
\end{mylemma}
\prooflem
Note that in view of scaling  (\ref{scale}), the supremum in formula
(\ref{kappa}) can be bounded
over all $x, y$ such that $\np{x-y}=1$ and. therefore, the denominator in
(\ref{kappa})
is greater than or equal to $ \Delta^{(m)}_s+ \Delta^{(\infty)}_s$.
On the other hand, since $\forall x \in \Qp$
\begin{eqnarray}
\abs{ \chi_{n}^{(\infty)}(x) - \chi_{n}^{(m)}(x) }=
\abs{1-\exp (\frac{ i2 \pi }{p} \sum_{k=m+1}^{\infty} x_{n-k} p^{-k}) }
< 2\pi p^{-m}
\end{eqnarray}
it is easy to show that the numerator in (\ref{kappa})  is always less than
$4\pi p^{-m}/(1-\abs{s})$.
\Endproof
Thus. Theorem \ref{mueqxi} permits the integration technique in $\Qp$
\cite{VVZ} to be applied for calculating integrals with respect to the
fractal measure $\mu_f$ or $\mu_f^0$.
For instance, it is possible to calculate the integrals of arbitrary
polynomials in  $z$~and~$\bar z$  with respect to $\mu_f^0$, for $z \in \sC$.
Indeed, it can be proved
\begin{eqnarray}\label{moment}
<z^L \bar z^{\bar L}> \equiv  \int \limits_{\sC} z^L \bar z^{\bar L}
 d\mu_f^0(x) =
\int \limits_{\Zp} (\Imm(x))^L( \overline{ \Imm(x)})^{\bar L}d\chi(x) =
\sum_{n,\bar n =0}^{\infty} C_{n,\bar n}^{L,\bar L} s^n {\bar s}^{\bar n}.
\end{eqnarray}
holds. Furthermore, $ C_{n,\bar n}^{L,\bar L} \in \Ns$ for $m=0,\infty$.
In particular, for $m=\infty$, the expression $C_{n,\bar n}^{L,\bar L}$
is the number of representations $n=n_0+...+n_L$~
and $\bar n=\bar n_0 +...+\bar n_{\bar L}$ such that
$ \sum_{k=0}^{L} p^{-n_k} =\sum_{k=0}^{\bar L} p^{- \bar n_k} {\rm mod}(p)$.
Applying this formula, we can show that, say,
$<z^p>=<\bar z^p>=1$
and $ <z^L \bar z^{\bar L}> =(1-\abs{s}^2)^{-L}\delta_{L,\bar L}$  for $L,
\bar L <p$.

Laying aside the problem of strict mathematical justification, we indicate one
more possible application
of Theorem \ref{mueqxi}. Consider a quantum particle on the fractal $ {\bf F
} =\Ran (\Imm)$.
In this case, $L^2(\sC,\mu_f)$ can be regarded as a Hilbert space of quantum
state.
Assume that the Schrodinger equation for energy eigenvalues $E$ can be
written in the form
\footnote
{
This quantum system can be thought of as an electron moving in a plane and
situated in a deep potential well supported by a fractal $F$. Then, in the
strong coupling limit, it can be assumed that the wave function is completely
localized on $F$ and the kernel $K$ is determined by overlap-type integrals
(see, e.g., \cite{Zim}).
}.
\begin{eqnarray}\label{DSEQ}
E \Psi(z) =\int K(\abs{z-\acute z}) \Psi(\acute z) d\mu_f(\acute z),
\end{eqnarray}
Using Theorem \ref{mueqxi} , we can pass to an equivalent equation in
$L^2(\Qp,\chi)$
with the kernel$K(\abs{\Imm(x)-\Imm(y)})$ for $x,y \in \Qp$. Evidently,
in the general case, this does not give any new results because the
kernel $K$ is translation-invariant in $\sC$, whereas $\mu_f$ is not.
Conversely, $\mu_f$ is invariant in $\Qp$, whereas $K$ is not.
However. applying  \bref{scale}, we can write the kernel in the form
\begin{eqnarray}\label{Kapp}
K( \npa{x-y}\cdot \abs{ ( \chi^{(m)}_{v(x-y)}(x-y)-1)  -O^{x,y}(s)}),
\end{eqnarray}

where $\abs{O^{x,y}(s)} \leq \abs{s}/(1-\abs{s})~ \forall x,y \in \Qp$.
Hence, for$\abs{s} \ll 1$ the kernel depends solely on $x-y$ and,
therefore, the Fourier transformation for $\Qp$ \cite{VVZ,GG} brings the
Hamiltonian to a diagonal form. In addition for $p = 2,3$,
the relation  $\abs{\chi^{(m)}_{v(x)}(x)-1 }=2\sin(\pi/p)$ can be applied
to prove that the
spectrum of the Hamiltonian (in this approximation) has the form
$\{ E_n=\widetilde K(p^{n}): n \in \Zs\} $and all eigenfunctions. except
for the ground state, are strictly localized
(i.e., they are compactly supported \cite{VVZ}). Here $\widetilde K(\cdot) $
is the Fourier image of the function $K(2\sin(\pi/p)\npa{\cdot})$.
Applying formulas(\ref{hina},\ref{impow}) and the $ p$-adic integration
technique, and consecutively expanding \bref{Kapp} into power series in
$s$ (but for fixed $a$), we can find the Fourier images of the coefficients
of the series at least in the case $m=\infty$ and $p=2,3$.
Thus, there appears to be a possibility of calculating subsequent corrections
to the spectrum and wave functions using perturbation theory. Note that the
foregoing equally applies to the diffusion problem on the fractal ${\bf F}$.
\section
{Continuous mappings of $\Tp$ into $\Rs^3$
}\label{immerstp}

For arbitrary $\forall n \in \Zs $~ and ~$\forall m \in \eNs \equiv  \Ns
\cup \{\infty\} $~
, we define some complex-valued continuous functions $\hins(\cdot)$  on$\Rs
\times \Zp$
by setting  $\forall (\xi,x)  \in \Rs \times \Zp$
\begin{eqnarray} \label{hinsol}
\hins(\xi,x)=\exp (i 2 \pi \xi \frac{\theta(p^{m-n}-\np{x})}{p^{\min(n,m)+1}})
\hin{x}.
\end{eqnarray}
Here the Heaviside function $ \theta(\cdot)$ is the indicator of the set
$[0,\infty]$.
It is easy to show that the functions  $\hins(\cdot)$ satisfy the relation
\begin{eqnarray}\label{xiconst}
\hins(\xi+l,x-l)=\hins(\xi,x),
\end{eqnarray}
$ \forall l \in \Zs, (\xi,x) \in \Rs \times \Zp$,
i.e., they are constant on the cosets of group $B$ (see Definition
\ref{soldef}). Note that
$\forall f =(\xi,x)  \in \Tp ~ \hinsn{\infty}(\xi,x)=\tilde
\chi_{(-p^{-n-1)})}(f)$,
where $\tilde \chi_{(q)}(\cdot) $
are continuous additive characters on $\Tp$ \cite{HR} (in view of Remark 4).
\begin{mydef}\label{mapTC}
For every $ s \in \UC $ and ~$\forall m \in \eNs $,~  we define the
continuous mapping
$ \wsol : \Rs \times \Zp  \mapsto \sC$,
by the formula
\begin{eqnarray} \label{wsolpow}
\wsol (\xi,x)  =\sum_{n=0}^{\infty} s^{n} \hins(\xi,x).~~~\forall (\xi,x)
\in \Rs \times \Zp
\end{eqnarray}
\end{mydef}
Note that  $\forall x \in \Zp $ the relation
\begin{eqnarray}\label{w0sz}
\Imm(x)=\wsol(0,x).
\end{eqnarray}
holds. For given $a \in \sC$, we introduce
the mapping $\sigma_a : \Su \times \sC \mapsto \Rs^3 $ by defining
the correspondence $\Su \times \sC \ni (e^{i2\pi \xi},z) \mapsto (x_1,x_2,x_3)
\in \Rs^3 $
according to the formulas
\begin{eqnarray}\label{cylmap}
x_1+i x_3 = e^{i 2 \pi \xi} \abs{a}  \left (1+ \RE{  \frac{z}{a} }   \right ),
x_2 =   \abs{a} \IM{\frac{z}{a}},
\end{eqnarray}
where $\Su$ is the unit circle in  $\sC$. We now define a mapping from  $\Tp$
into  $\Rs^3$
as the composition of $\sigma_a  $  a mapping from $\Tp$ ¢ $\Su \times \sC $
is constructed below.
To this end, we consider $\Tp$ and $\Su \times \sC$ as quotient spaces with
respect to the action of the group $\Zs$ on$\Rs \times \Zp $~ and~ $\Rs \times
\sC$,
namely, $\Rs \times \Zp $~ ¨~ $\Rs \times \sC$,
the action on  $\Rs \times \Zp $~  and  ~ $\Rs \times \sC$ is defined as
shifts by the elements
$(n,-n) \in B$ and  $(n,0) \in \Zs\times \{0\}$, respectively. Let
$\idwsol \equiv id \times \wsol :\Rs \times \Zp \mapsto \Rs \times \sC $,
where $id$  is the identity mapping.
It is easy to see that property \bref{xiconst}
guarantees the existence of the continuous quotient mapping $\facwsol: \Tp
\mapsto \Su \times \sC $, as well as a mapping $\Wsol: \Tp \mapsto \Rs^3$,
such that the diagram
\\
\begin{picture}(120,130)(-100,-110)
\put(63,-0) {\makebox(0,0){ $\idwsol$ }}
\put(25,-10) {\makebox(0,0){$\Rs \times \Zp$ }}
\put(110,-10) {\makebox(0,0){$\Rs \times \sC$ }}
\put(45,-10){\vector(1,0){35}}

\put(27,-20){\vector(0,-1){22}}
\put(102,-20){\vector(0,-1){22}}
\put(20,-30){\makebox(0,0) {$\phi$}}
\put(110,-30){\makebox(0,0) {$\varphi$}}

\put(63,-40) {\makebox(0,0){ $ \facwsol $   }}
\put(30,-50) {\makebox(0,0){$\Tp$ }}
\put(110,-50) {\makebox(0,0){${\bf S^1} \times\sC$ }}
\put(40,-50){\vector(1,0){40}}
\put(102,-60){\vector(0,-1){30}}
\put(40,-65){\vector(3,-2){40}}
\put(50,-85) {\makebox(0,0){$\Wsol$}}
\put(110,-75) {\makebox(0,0){$ \sigma_a$}}
\put(105,-100) {\makebox(0,0){$\Rs^{3}$ }}
\end{picture}
\\
is commutative. Here $\phi$~and~$\varphi$ are the corresponding
canonical projections.
Moreover, it can be shown that the actions of the group
$\Zs$ on $\Rs \times \Zp $~ and ~ $\Rs \times \sC $
are compatible in the sense that $\forall f \in \Tp$ the relation
\begin{eqnarray}\label{accord}
\idwsol( \phi^{-1}(f)) = \varphi^{-1}(\facwsol(f) ) .
\end{eqnarray}
holds.
We introduce an invariant metric $\hat d$ on $ \Rs \times \sC $
such that $$\hat d( (\xi_1,z_1),(\xi_2,z_2)) =
\max (\abs{\xi_1-\xi_2},\abs{z_1-z_2} )
~ \forall (\xi_1,z_1),(\xi_2,z_2) \in  \Rs \times \sC$$
The metric $d$ on $ \Su \times \sC $ is defined as the quotient metric of
$\hat d$ following the same scheme as in \bref{rt}.
Then it is not difficult to prove that $~ \forall
(e^{i2\pi \xi_1},z_1),(e^{i2\pi \xi_2},z_2) \in  \Su \times \sC$
(here we assume that $\xi_1,\xi_2 \in [0,1)$),
\begin{eqnarray}
d( ( e^{i2\pi \xi_1},z_1),( e^{i2\pi \xi_1},z_2)) =
\max ( \inf_{n=0,\pm 1 } (\abs{\xi_1-\xi_2+n}),\abs{z_1-z_2} ).
\end{eqnarray}
Because $\phi,\varphi$ and the restriction of $\sigma_a$
to any bounded set in $ \Su \times \sC$  are L-contractions,
if we show that $\idwsol$ is an L-contraction. then this implies that
$\facwsol$ and
$\Wsol$  are L-contractions as well.
Moreover. if $\idwsol$ is an L-isometry. then so is $\facwsol$ . Using
\bref{accord},
we can show that $\forall f,g \in \Tp$, there are $(\xi,x) \in
f,~(\eta,y)\in g$  such
that, $\forall \ep >0$, the chain of inequalities
\begin{eqnarray*}
\rho_{\alpha}(f,g) \leq \hat \rho_{\alpha}((\xi,x) ,(\eta,y))
\Ll  \hat d(\idwsol(\xi,x),\idwsol(\eta,y))
\leq d(\facwsol(f),\facwsol(g))+\ep .
\end{eqnarray*}
holds. It is also obvious that if the restriction of $\sigma_a$ to
$\Ran( \varphi \circ \idwsol)~$ an L-isometry, then the mapping
$\Wsol$ is also.
To find the conditions under which $\idwsol$ is an L-contraction or an
L-isometry, we introduce the following two numbers:
\begin{eqnarray}
\tilde \Delta^{m}_s=\inf_{x,y\in \Zp, \xi \in \Rs}
(    \abs{ s^{v(x-y)} (\wsol(\xi,x)-\wsol(\xi,y))} ),\\
\gamma^{m}_{s,a}= - \inf_{(\xi,x) \in \Rs \times \Zp} (
{\rm Re}(\wsol(\xi,x)/a)).
\end{eqnarray}
It is easy to show that $\tilde \Delta^{m}_s$  satisfies the same
inequality as the one in \bref{ds0}
and, hence, $\tilde \Delta^{m}_s>0$ for $\abs{s} < s_0$.
Furthermore, it can be proved that $~\tilde \Delta^{\infty}_s=
\Delta^{\infty}_s$ (see below).
The above-mentioned conditions follow from the chain of inequalities
below:
\begin{eqnarray*}
\tilde \Delta^{m}_s \max ( \abs{\xi-\eta},\np{x-y}^{\frac{1}{D_s}} ) \Ll
\max ( \abs{\xi-\eta}, \Delta^{m}_s \np{x-y}^{\frac{1}{D_s}}-
\frac{2\pi}{p(1-\abs{s})} \abs{\xi-\eta} )\\
\leq
\max ( \abs{\xi-\eta},\abs{\wsol(\xi,x)-\wsol(\eta,y) }  ) \\ \leq
\max ( \abs{\xi-\eta},\frac{1}{(1-\abs{s})}  ( \np{x-y}^{\frac{1}{D_s}}+
\frac{2\pi}{p} \abs{\xi-\eta}) )
\Ll
\max ( \abs{\xi-\eta},\np{x-y}^{\frac{1}{D_s}} ) ,
\end{eqnarray*}
the first and last being the consequences of relations \bref{LEQ}; the
others were obtained analogously to \bref{AA}. An elementary geometric
consideration implies that the restriction of aa to any compact set in
$\Pi_a \equiv \{z \in \sC : \RE{z/a} > 1 \}$ is an L-isometry and,
furthermore, $\Ran( \facwsol) \subset \Pi_a$ for $\gamma^{m}_{s,a} < 1$.
Thus, we have proved the following theorem.
\begin{mytheorem}\label{LES}
$\hat \rho_{\alpha},~ \rho_{\alpha},~\hat d$ and $d$
be the metrics in
$\Rs \times \Zp,\Tp, \Rs \times \sC $ and $\Su \times \sC$, respectively.
Let $\Rs^3$ be endowed with the standard Euclidean metric and let
$\alpha =D_{s}^{-1}$.
Then for any $ m \in \eNs $  and $s \in \UC$. the following assertions hold:
\\$1)$ $\idwsol$  is an L-contraction and, therefore, $\facwsol$~and~$\Wsol$
are also L-contractions:
\\$2)$ if $\tilde \Delta^{m}_s > 0$, then $\idwsol$ and $\facwsol$ are
L-isometries
and  addition $\gamma^{m}_{s,a} < 1$, then $\Wsol$ is an L-isometry.
\end{mytheorem}
It follows immediately that if
{\it $s$ and $a$ are such that $\tilde \Delta^{m}_s > 0$ and $~\gamma^{m}_{s,a}
< 1$,
then the mapping $\Wsol$
a continuous embedding of $\Tp$ in $\Rs^3$
and the local and global Hausdorff dimensions $\Wsol(\Tp)$  are equal to
$D_s+1$.}

Next, we describe the geometric structure of the set   $\Wsol(\Tp)$.
To this end, we note that $(\Tp,\rho_{\alpha})$ can be regarded as the
total space of the locally trivial fiber bundle $(\Tp,\Su)$~ with the
projection $\chi_0: \Tp \mapsto \Su $, such that $ \chi_0 (\xi,x) =
\exp(i2\pi\xi) \in \Su,~\forall (\xi,x)\in \Tp $
and fibers $\chi_{0}^{-1}( \exp(i2\pi\xi))=(\xi,\Zp)$ isometric to
$(\Zp,\np{\cdot}^{\alpha})$.
It is easily seen that the mapping $\facwsol$ is a fiber morphism of the
bundle $(\Tp,\Su)$
into the trivial bundle $( \Su \times {\rm D}(r_s),\Su) \subset (
\Su \times \sC,\Su)$
(where ${\rm D}(r_s)$ is a closed disk of radius $r_s=(1-\abs{s})^{-1}$  ).
In this case. the fibers are the sets $\wsol(\xi,\Zp) \subset {\rm D}(r_s)$
for every fixed $\xi$ and if $\facwsol$  is an L-isometry, then these sets
are fractals with global and local dimensions equal to $D_s$.
Moreover, it follows from \bref{w0sz} that the fiber with  $\xi=0$ is
isometric to $\Imm(\Zp)$.
For simplicity, we assume that $a > r_s$, Then  $\sigma_a$ is the standard
embedding of $( \Su \times {\rm D}(r_s),\Su)$ in a solid torus   in
$\Rs^3$ (a torus together with its interior) and the fibers are mapped
isometrically onto disks lying in the plane turned around the $x_2$ -axis
through an angle $2\pi \xi$ relative to the plane  $x_3=0$.

Let us define the action of an element $t$ of the group  $\Rs$  on
$f \in \Tp$  as the shift  $f  \mapsto f_t = f+\phi(t,0)$ and consider
the orbits of this group in $\Tp$.  It can be shown \cite{HR} that each
orbit is a dense subset in $\Tp$.
Because the mappings  $\facwsol$~and~$\Wsol$ are L-contractions, the
images of these orbits are continuous curves that also wind densely
around  $\facwsol(\Tp) \subset  \Su \times {\rm D}(r_s) $~
and   ~$\Wsol \subset \Rs^3 $ , respectively.

It turns out that in the case $m=\infty$ the images of the orbits of the
group $\Rs$ are smooth curves that are trajectories of a dynamic system.
\begin{mytheorem}
Let $\Wsoli$  be an L-isometry. Then there is a global
Lipschitz vector field $\Gamma: \Rs^3 \mapsto \Rs^3$ such that
$\forall f \in \Tp$,
\begin{eqnarray}\label{gamma}
\Gamma( \Wsoli( f))= -2\pi \left (
L_{(2)} \Wsoli( f)  + p^{-1}{\rm Re} \left(\tilde
\chi_{-1}^{(\infty)}(f)( L_{(3)} + i L_{(1)}) \right) \Wsolip(f) \right ).
\end{eqnarray}
The action of the one-parameter homeomorphism group  $U_t: \Rs^3 \mapsto
\Rs^3$, generated by the equation
\begin{displaymath}
{\rm \frac{d}{dt}}{\bf r} = \Gamma({\bf r}),
\end{displaymath}
is compatible with the action of $\Rs$ in $\Tp$ in the sense that
$\forall f \in \Tp, \forall t \in \Rs$,
\begin{eqnarray}\label{Hom}
\Wsoli( f_t)=U_t \Wsoli(f).
\end{eqnarray}
Here  $(L_{(k)})_{i,j}=\ep_{i,j,k}$, where $\ep_{i,j,k}$ is the Levi-Civita
symbol and
$ \Wsolip(f) \equiv \lim \limits_{\ep \rightarrow  +0}  \Wsolipa(f)$.
\end{mytheorem}
\proof
First, we show that $\Gamma_{|\Ran(\Wsoli)}$ is an L-contraction. Obviously,
it is sufficient to prove that
$((\Wsoli)^{-1} \circ \Wsolip)_{|\Ran(\Wsoli)}$ is an L-
-contraction.
However, this is a direct consequence of the fact that $(\Wsoli)^{-1}$ ¨
$ \Wsolip$
are L-contractions for $( \Tp , \rho_{\alpha} )$ if $\alpha = D_s^{-1}$
(  $ \Wsolip$   is an L-contraction since $\rho_{\alpha+1} \Ll \rho_{\alpha}$).
The existence of a Lipschitz mapping $\Gamma$
defined throughout $\Rs^3$ and satisfying  \bref{gamma}
is now
ensured by the Kirszbraun theorem \cite{FED}.
It follows immediately from formulas  (\ref{hinsol},\ref{wsolpow}) that
$\forall t \in \Rs ,x \in \Zp$ the equation
\begin{eqnarray}
\frac{d}{dt} \wsoli(t,x)=\frac{i2\pi}{p}\wsolip(t,x).
\end{eqnarray}
hold.  With the help of  \bref{cylmap}, it can be deduced from this equation
that $\forall f \in \Tp$
and $\forall t \in \Rs$ we have
\begin{displaymath}
\frac{d}{dt}  \Wsoli  ( f_t)=\Gamma ( \Wsoli(f_t) ).
\end{displaymath}
According to theory of differential equations \cite{ARN}, this implies formula
\bref{Hom}.
\Endproof
Let us mention another specific peculiarity of the case $m=\infty$.
It follows from formulas (\ref{hina},\ref{impow},\ref{hinsol},\ref{wsolpow})
that $\forall x \in \Qp$,
\begin{eqnarray}\label{vqwt}
[\Immi](x)=\wsoli(\fp{x},\ip{x}).
\end{eqnarray}

This implies the above-mentioned relation $~\tilde \Delta^{\infty}_s=
\Delta^{\infty}_s$.
Moreover, using Theorem \ref{JQT} and the Kirszbraun theorem \cite{FED},
we can now prove that
if $\Wsoli$ is an L-isometry, then there exists an L-contraction $J :
\sC \mapsto \Rs^3$
such that $J$ isometrically maps every cluster $B_{l}^0=\Immi(\{ x \in \Qp :
\np{x-l}<1 \}$
onto the image of the $\fp{l}$-fiber  $\Wsoli( (\fp{l},\Zp) )$  and we have
\begin{eqnarray}\label{slice}
\Wsoli \circ j =J \circ \Immi.
\end{eqnarray}

\section{   Explanation of the figures}\label{pic}

\begin{figure}
\centerline{\epsfxsize=14 truecm \epsfbox{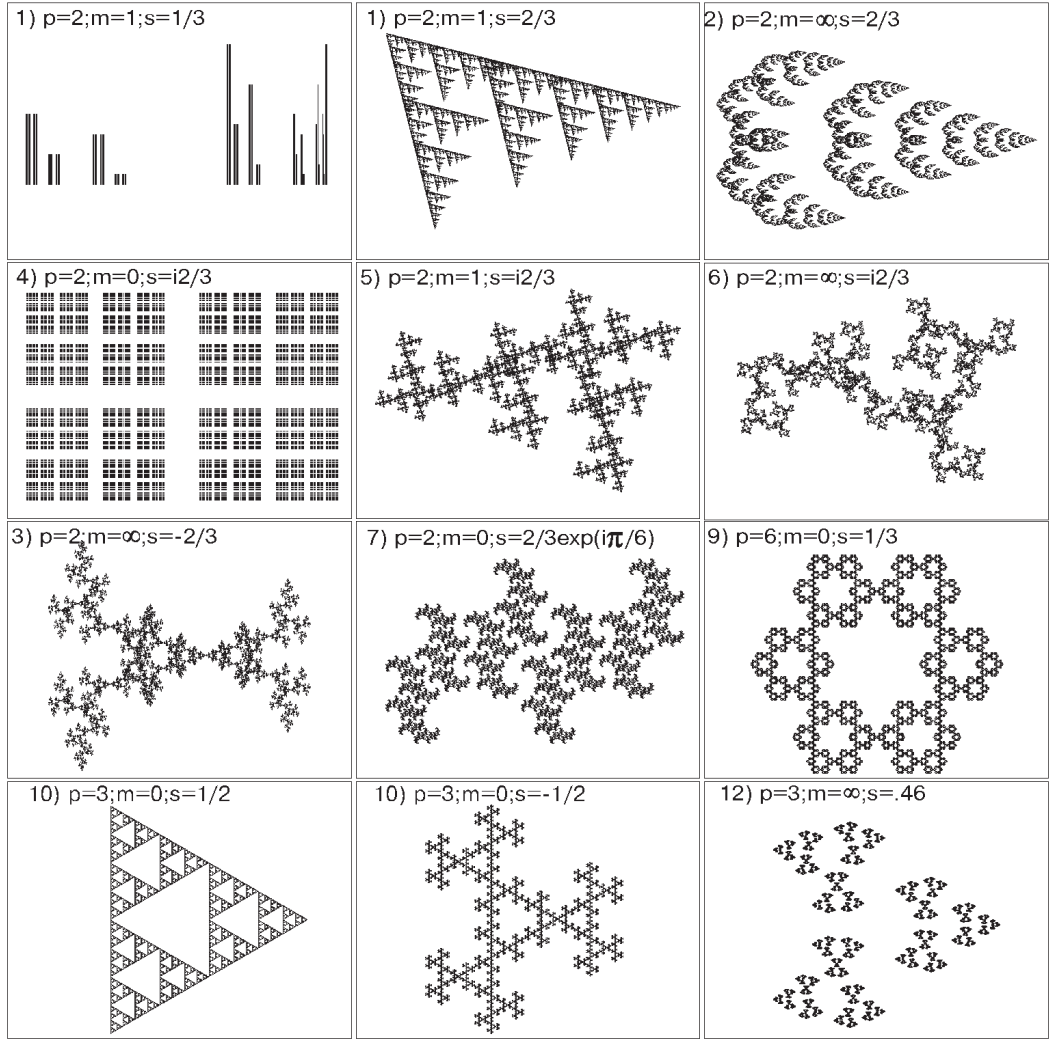}}
\caption{  Images of $\Zp$ in the complex plane for different values
$p,m$ and $s$}\label{fpi}.
\end{figure}
\begin{figure}
\centerline{\epsfxsize=14 truecm \epsfbox{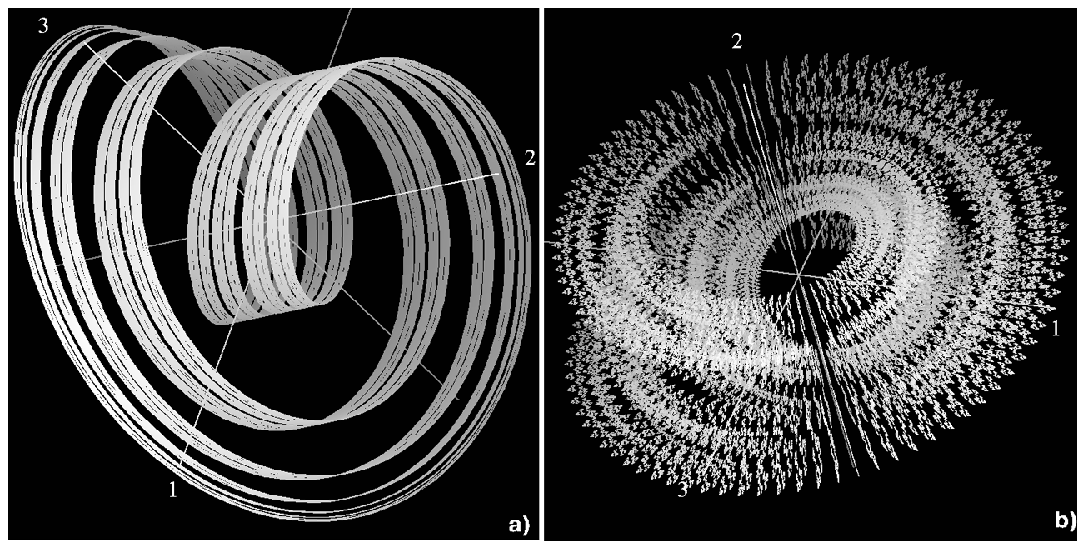}}
\caption{   Embeddings in $\Rs^3$~$a)$  $\Tn{2},~(m=0)$;
$b)$~$\Tn{3}$ (by fibers) or $3^{-4}\Zn{3}$ ,($m= \infty$).}\label{spi}
\end{figure}

Due to the rapid convergence of series  \bref{power}  and  \bref{impow}
and, also, because even not very large positive integers form a sufficiently
dense network in $\Zp$, the computer construction of the sets
~$\Imm(\Zp)$~and $\Wsol(\Tp)$ encounters no fundamental difficulties.
Figure 1.1 represents a Cantor set and. for the sake of clarity, a segment,
whose length is proportional to one of the numbers $e_i =1,3,5,7,2,6,10,14$
is associated with each point if the point of the image has the form
$\Immnn{0}{1/3}(e_i \cdot y^2)$ for  $y \in \Zp$ (for these numbers, see
\cite{VVZ}).
The set in Fig. 1.10 is obviously the Sierpiriski  triangle. It can be
seen from Fig. 1.9 that the boundaries of the connected components of the set
$\sC \backslash \Immnn{0}{1/3}(\Zn{6})$
consist of Koch curves (note that $p$ is not a prime number: see footnote 3).
Finally, Fig. 1.4 can serve as an illustration of the fact that $\Zn{4}$ is
homeomorphic to $\Zn{2}$.
Figure 2a demonstrates an embedding of $\Tn{2}$  into $\Rs^3$ with the
parameters $s=1/2.2,~ a=i2$ and $m=0$. and the Cantor structure of this set
can be distinctly seen.
Figure 2b illustrates an embedding
of $\Tn{3}$ ($s=s_0 - 0.02 \approx 0.46 ,~ a=5/2$ and $m=\infty$ ).
Fibers for which the values of  $\xi$ are multiples of $1/81$ are represented
and, therefore, Fig. 2b can simultaneously be regarded as the image of
$\Immi(3^{-4}\Zn{3})$ under the mapping J (see formulas  (\ref{vqwt},
\ref{slice})); this image is, in fact, shown in Fig. 1.12.
In conclusion, we note that additional general constructions, such as a-adic
numbers and solenoids (see [11]), can also be treated in a similar way, at
least for m = $\infty$. It seems, in this case, that constructing an embedding
of ${\bf A} / \Qs$ in $\Rs^3$, for instance (where ${\bf A}$  is the adele ring
\cite{Mam}), could be interesting.

The research was partially supported by the RFBR
grant 01-02-17682-a and by the INTAS grant 00-00334.


\begin{thebibliography}{99}
\bibitem{Kob}
{\it N. Koblitz}\ p-Adic Numbers, p-Adic Analysis, and Zeta-Functions,
Springer, New York Heidelberg Berlin (1977).
\bibitem{GG}
{\it I. M Gel'faiid, M. I. Graev. and L. I. Pyatetskii-Shapiro }\
Representation Theory and Automorpttic Functions. Saunders. Philadelphia
(1969).
\bibitem{VVZ}
{\it  V. S. Vladimirov, I. V. Volovich, and E. I. Zeienov}
\ p-Adic Analysis and Mathematical Physics, World Scientific, Singapore-New
Jersey-London-Hong Kong (1994).
\bibitem{FED}
{\it H. Federer }\
Geometric Measure Theory, Springer, New York-Heidelberg-Berlin (1969).
\bibitem{HR}
{\it E. Hewitt and K. Ross}\
Abstract Harmonic Analysis, Vol. I, Springer, New York-Heidelberg Berlin (1
\bibitem{bill}
{\it P. Billingsley} \ Ergodic Theory and Information, Wiley, New
York-London-Sidney (1965)
\bibitem{Feder}
{\it J. Feder} \ Fractals, Plenum, New York (1988).
15. , \bibitem{Mam}
{\it D. Mumford.} \ Tata Lectures Notes on Theta Functions, Vols. I, II,
Birkhauser, Boston-Basel-Stuttgart (1983. 1984).
\bibitem{Kel}
{\it J. L. Kelley} \ General Topology. Van Nostrand, Princeton,
New Jersey (1957).
\bibitem{Shabat}
{\it B. V. Shabat} \ Complex Analysis [in Russian], Vol. 1, Nauka,
Moscow (1985).
\bibitem{Zim}
{\it J. M. Ziman} \ Models of Disorder. The Theoretical Physics of
Homogeneously Disordered Systems. Cambridge University Press,
Cambridge-London-New York-Melbourne (1979).
\bibitem{LL}
{\it A. J. Lichtenberg and M. A. Lieberman} \ Regular and Stochastic
Motion, Springer, New York- Heidelberg Berlin (1983).
\bibitem{ARN}
{\it V. I. Arnold} \ Ordinary Differential Equations, MIT Press,
Cambridge (1978).
Œ.:  ãª , 1975.
\bibitem{FFT}
{\it K. Alien and M. Kluater}\ Optical Transformations in Fractals.
Fractals in Physics, North Holland. Amsterdam Oxford-New York-Tokyo (1986).
\bibitem{Zelenov}
{\it Zelenov E.I. }//
J.Math.Phys. V32.147-152. 1991.
\bibitem{Pit}
{\it Pitkanen M.} // $p$-adic Physics ?. Department of Theoretical Physics,
University of Helsinki, SF-00170 Helsinki, Finland. 8. September 1994.
\bibitem{dif1}
{\it Havlin S., Weissman N. }//
J.Phys.A 19. L1021-1026.1986.
\bibitem{dif2}
{\it Ogielski A.T., Stein D.L. }//
Phys. Rev. Let. V.55. N15. 1985.
\end{thebibliography}
\end{document}